\documentclass{article}
\usepackage{amssymb,amsmath, tikz}
\usepackage{graphics}
\input xy
\xyoption{all}
\usepackage{braket}
\usepackage{hyperref}
\newlength{\dhatheight}

\newcommand\phantomarrow[2]{%
  \setbox0=\hbox{$\displaystyle #1\to$}%
  \hbox to \wd0{%
    $#2\mapstochar
     \cleaders\hbox{$\mkern-1mu\relbar\mkern-3mu$}\hfill
     \mkern-7mu\rightarrow$}%
  \,}

\begin{document}

\title{Separation and symmetry on two dimensional manifolds}

\author{\v{S}t\v{e}p\'an Hude\v{cek}, Svatopluk Kr\'ysl\footnote{{\it E-mail address}: 
Svatopluk.Krysl@mff.cuni.cz}\\ {\it \small Charles University, Faculty of 
Mathematics and  Physics,  Czechia}
\thanks{The second author thanks for financial supports from the founding 
No. 20-11473S  granted by the Czech Science Foundation.}}

\maketitle \noindent
\centerline{\large\bf Abstract} 

We introduce notions of a separated solution and of a simple symmetry that generates a   differential 
operator on a smooth manifold. 
We prove that a  differential operator on a two dimensional manifold has a separated solution if it has 
a homogeneous simple symmetry of degree one which does not generate the operator.

{\it Keywords:} Symmetry operators; separation of variables; separation on manifolds; Bessel functions

{\it Classification:} 58J70, 35A30, 94A11

\section{Introduction}

We investigate separable solutions of partial differential operators on smooth manifolds. Our aim is a precise formulation and a proof of   a separation theorem.
It is well  known that the problem of  variables' separations  is connected to symmetries.
See  Miller \cite{Miller}, Kostant \cite{Kostant} and Eastwood \cite{East} for varied approaches.
Let us remark that the notion of a separable operator is not precisely defined in the book of Miller \cite{Miller}, although  
assertions concerning this notion are contained there. In Koornwinder \cite{Koorn}, we find a broad  discussion on what   an appropriate definition of a separable operator shall be and a 
history of this problem. According to \cite{Koorn}, the  history of the `separation problem' goes back to St\"ackel \cite{Stackel}. 
Koornwinder writes in the ibid. citation that he is motivated by 
a missing definition of a separable operator in Miller \cite{Miller}. This was one of the motivations for authors of this paper as well. 
Miller gives a definition of a separable operator, which is subsequently, used and investigated. See, e.g., Hainzl \cite{Hainz}.

  In particular we are interested in smooth manifolds of dimension two.  Let us remark that we define neither a separable differential operator, nor a separable equation. 
	We   define the notion of a  {\it separable solution} with respect to a chosen coordinate system instead.
Our definition of a {\it symmetry operator} of a differential operator is slightly different from the definition in Miller \cite{Miller}, who set it in the case of the Helmholtz and Laplace operators. Namely, we do not bound the order of the symmetry operator.
The definition is also different  from  the definitions of Eastwood \cite{East} and Kostant \cite{Kostant}.
These definitions are compared in the paper briefly. 

We also  define the notions of the  {\it invariance of an operator with respect to a smooth function}, and give conditions under which a symmetry operator is said to generate a given differential operator.   In the case of the Helmholtz operator on the Euclidean spaces $\mathbb{R}^2$ and 
$\mathbb{R}^3,$  Miller  in \cite{Miller} does not treat the problem explicitly.  
We prove that a differential operator $D$ defined on a two dimensional manifold  
has a {\it non-constant} separated solution 
if there is a 
homogeneous degree one symmetry  which is so called  simple  and     does not generate $D$ (Theorem 3). Let us remark that
a homogeneous degree one symmetry is a vector field. If we omitted in the formulation of the theorem that the solution is non-constant, the assertion would be almost trivial because the zero solution is always separated.  
At the end of the paper, we apply our results in the well known case of the Helmholtz operator on the Euclidean plane.

\section{Separation theorem}

Let $M$ be a smooth manifold of   dimension $n$ and let us denote the filtered   ring of complex {\it linear} differential operators on $M$ by $DO(M).$ The filtration is given by the order of the operator.  See, e.g., Seeley \cite{Seeley} for a definition of a differential operator on a manifold.  The ring $DO(M)$ is non-commutative
if $n\geq 1.$ The word linear is usually omitted in differential geometry. We follow this convention.
The multiplication in the ring is the composition  and the addition is the   addition of operators.
The ring $DO(M)$ is a left module over  the ring $C^{\infty}(M)$ of smooth complex valued  functions on $M.$ 
The module structure is defined by setting $(f\cdot D)(g) = fD(g)$ for $f, g \in C^{\infty}(M)$ and $D \in DO(M),$ where at the right 
hand side the point-wise multiplication of functions is considered. In the paper, we omit writing the dot in $f \cdot D.$
   
 Let $R$ denote the subring of the ring $C^{\infty}(M)$ consisting of constant functions. The ring $DO(M)$ is also an {\it $R$-algebra}, i.e., 
it is a ring and an $R$-module 
satisfying the compatibility condition
$r(D_1D_2)=(rD_1)D_2=D_1(rD_2)$ for each $D_1, D_2 \in DO(M)$ and $r \in 
R.$ We identify $R$ with the ring $\mathbb{C}$ of complex numbers.

\vspace{0.5cm}

{\bf Notation:} Let us recall that an ordered $n$-tuple with entries in
$\mathbb{N}_0=\mathbb{N} \cup \{0\}$ is called a multiindex of length $n.$ For a multiindex $I=(i_1,\ldots, i_n),$ we set
$|I|=\sum_{j=1}^n i_j.$  We consider the component-wise addition and subtraction of multiindices.

\vspace{0.5cm}

If $(U,\phi)$ is a map in the   atlas of the manifold $M$, we define functions $x^i: U\to 
\mathbb{R},$ $i=1, \ldots, n$ (the coordinate functions) by 
$\phi(m)=(x^1(m),\ldots, x^n(m)),$ where $m \in U.$
We write $\partial_{x^i}$ for the appropriate local vector field $\frac{\partial}{\partial x^i}$ which is also 
a first order   differential operator. For a non-negative integer $a,$ $\partial_{x^i}^a$ denotes the $a$-fold composition of $\partial_{x^i}.$ (If $a=0,$
we take $\textrm{Id}_{C^{\infty}(M)}.$)
If $I=(i_1,\ldots, i_n),$  we use the symbol 
$\partial_I$ for $\partial_{x^1}^{i_1} \ldots \partial_{x^n}^{i_n},$ where the right-hand side is  
  the multiplication  of the appropriate differential operators. We assume that the order of positions of $\partial_{x^j}^{i_j}$ in $\partial_I$ is fixed by the order of the coordinate functions of the chosen  map. (Thus, e.g.,
$\partial_{x^2}\partial_{x^1}$ cannot appear in the expression for $\partial_{I}$ although it equals to the allowed $\partial_{x^1}\partial_{x^2}.$)
In the case of $\mathbb{R}^2,$ we   denote the canonical 
Cartesian coordinates by $x^1$ and $x^2$ or by $x$ and $y$.

\vspace{0.5cm}

{\bf Remark:} Notice that when writing $D=\sum_I f^I \partial_I$ (finite sum), even locally on an open set $U$,  
the coefficients $f^I$ are uniquely defined by a choice of the manifold map. We write the functions first and then the composition of the vector fields. 
 
For a chosen map, we can  decompose $D$ of order $d$ as $\sum_{i=0}^{d}\sum_{|I|=i}f^I \partial_I$. For each $i,$
the sum $\sum_{|I|=i}f^I \partial_I$ is called the homogeneous component of $D$ of degree $i$.   A differential operator is called a {\it homogeneous operator of degree $i$} if it equals to a single homogeneous component of degree $i$ only. Note that the notion of the homogeneous operator  depends on coordinates (see the example of the Laplace operator below).
However, the zero degree component  is independent on coordinates, since it is given by the evaluation of the operator on the constant function $1.$
Independently on maps, we may  speak about {\it first order homogeneous operators} as those first order operators whose homogeneous component of degree zero vanishes. These are actually vector fields on the manifold.

We refer also to Bj\"ork \cite{Bjork} for algebraic properties of the subalgebra  of $DO(\mathbb{R}^n)$ which consists
of operators which are of the form $\sum_I p^I \partial_I$ (finite sum), where $p^I$ are polynomials with respect to the  Cartesian
coordinates on $\mathbb{R}^n,$ the so called Weyl algebra.

\vspace{0.5cm}

{\bf Definition 1:} Let $(U,\phi)$ be a map on $M$ and  $1\leq k \leq n.$ 
We say that $f \in C^{\infty}(U)$ {\it does not  
$\phi$-depend} on the first $k$ variables, if 
$$(f \circ \phi^{-1})(x_1',\ldots, x_k',x_{k+1},\ldots, x_n) = (f \circ 
\phi^{-1})(x_1'',\ldots, x_k'',x_{k+1},\ldots, x_n)$$ for any 
$(x_1',\ldots,\:x_k',\:x_{k+1},\ldots,\:x_n),\:(x_1'',\ldots,\:x_k'',\:x_{k+1},
\ldots,\:x_n) \in \phi(U).$ Similarly we define  a function that
{\it does not $\phi$-depend} on the last $k$ variables.
 
\vspace{0.5cm}

{\bf Remark:}
\begin{itemize}
\item[1] Degenerate cases. If  $f$ does not 
$\phi$-depend on the first $n$ or on the last $n$ variables, then $f$ is locally constant.
Naturally, we are interested in cases when $f$ is not locally 
constant, which leads to $k<n.$ 
\item[2.] If $f$ does not $\phi$-depend on the first $k$ variables, we have
$\partial_{x^1}f=\ldots =\partial_{x^k}f=0$ where $\phi=(x^1,\ldots, x^n).$ The opposite implication holds   locally.
Moreover, if $M \subseteq \mathbb{R}^n$ and $U$ is a convex set, the implication can be reversed.
\end{itemize}

\vspace{0.5cm}

{\bf Definition 2:} Let $(U,\phi)$ be a map on $M$ and $g$ be a smooth function on $U$ which 
does not $\phi$-depend on the first or on  the last  $k$ variables.
We say that a differential operator $D$ {\it is $g$-invariant} with respect to $(U,\phi)$ if there exists a 
differential operator $D'$ such that 
$D(gh)=gD'(h)$ for any smooth function $h$ on $U$ which does not $\phi$-depend on the
last or on the first $n-k$ variables, respectively.

\vspace{0.5cm}
 
{\bf Remark:}
 If $g$ is locally constant, then any differential operator is 
$g$-invariant (with respect to any manifold map) since one can take $k=n$ and $D'=D.$

\vspace{0.5cm}

{\bf Notation:} If $(U,\phi)$ is a map, $f: U \to \mathbb{C}$ is a function on 
$U$ and $(x^1_0, \ldots, x^n_0) \in \phi(U)$ is a point in $\mathbb{R}^n,$ we write $f(x^1_0,\ldots, x^n_0)$ briefly instead 
of 
$(f \circ \phi^{-1})(x^1_0,\ldots, x^n_0).$

\vspace{0.5cm}

{\bf Example 1:} Operator $D'$ from Definition 2 (the case of {\it first} $k$ variables) may   contain partial derivatives with respect to some of 
the last $n-k$ variables and   still fulfil the definition.
Namely, let $M = \mathbb{R}^2 \setminus (-\infty,0] \times \{0\}$ be equipped with the polar patch $(r,\theta) \mapsto (r \cos\theta,  r \sin\theta) \in M,$ where $(r,\theta) \in (0, \infty) \times (-\pi,\pi).$ We denote the appropriate inverse of the patch by $\phi$ and set 
$f(r, \theta)=e^{\lambda \theta},$ $\lambda \in \mathbb{C}.$ The function does not $\phi$-depend
on the first variable.
Since the Euclidean Laplace operator  $\Delta = \partial_{x}^2 + \partial_{y}^2,$ restricted to smooth functions on $M,$ equals   to 
$\partial_{r}^2 +\frac{1}{r} \partial_r + \frac{1}{r^2} \partial_{\theta}^2,$
 we see that for any function $h: r \mapsto h(r),$ $r \in (0, \infty),$ which does not $\phi$-depend on the last coordinate
\begin{equation*}
    (\partial_{r}^2 + \frac{1}{r} \partial_r+\frac{1}{r^2} \partial_{\theta}^2) (e^{\lambda \theta} h) =  
		e^{\lambda \theta} (\partial_r^2 + \frac{1}{r} \partial_r 
 + \frac{\lambda^2}{r^2})h
= e^{\lambda \theta} D'(h),
\end{equation*}
where $D' =  \partial^2 
+  \frac{\lambda^2}{r^2}.$ 
Thus for any $\lambda \in \mathbb{C},$ $\Delta$ is $e^{\lambda\theta}$-invariant.
We see that  $D'$ contains only partial derivatives in the variable $r.$ Nevertheless, we can
 take $D'= \partial_r^2 + \frac{1}{r} \partial_r  +  
\frac{\lambda^2}{r^2} + D'',$ where $D''$ is a differential operator containing   partial 
derivatives only with respect to $\theta$ and   whose zero order term is zero. In particular, we see that
 $D'$ may depend on the last variable $\theta,$ and also that it is not unique.

\vspace{0.5cm}

Any associative algebra $A$ over a field gives rise to the Lie algebra $(A, [,])$ with a Lie bracket that equals
to the {\it commutator}, which is   defined by $[C,D]=CD-DC$ for elements $C,D \in A.$

\vspace{0.5cm}

{\bf Definition 3:} Let $D$ be a differential operator. We call a differential 
operator $L$ a {\it symmetry operator} of $D$ if there exists a differential 
operator $G$  
such that $[D,L]=  GD.$  
We call $L$ a {\it simple symmetry operator} of $D$ if $[D,L]=0.$ 
We say that a differential operator $L$ {\it generates}  a 
differential operator $D$ if there exist an open set $U$ and smooth functions 
$f_0, \ldots, f_m$ defined on $U$
such that the restriction of $D$ to $C^{\infty}(U)$ 
satisfies $D =\sum_{i=0}^m f_i  L^i,$ where $L^0$ denotes
the identity operator on $C^{\infty}(U).$

\vspace{0.5cm}

{\bf Remark:} 
\begin{itemize}
\item[1.] The definition of a symmetry operator in Eastwood \cite{East} is more general than our definition.
The above definition of a symmetry operator generalizes the 
definition of Miller 
\cite{Miller}, pp. 2 and 7.  
The symmetry operator as defined in Kostant 
 \cite{Kostant}, p. 101, is less general than the symmetry operator  of Miller in \cite{Miller} if the formula in the Miller's definition
were considered for a general partial differential operator $D.$  

\item[2.] For an open set $U,$ we consider the set of differential operators $DO(U)$ as a 
left $C^{\infty}(U)$-module as explained above. For a differential operator $L$ 
let us consider the left $C^{\infty}(M)$-hull  (finite sums) of 
the the set $\{L^i, i \in \mathbb{N} \cup \{0\}\}$  and denote it by $\langle L \rangle.$
Then ``$L$ generates $D$'' means precisely  that  there is an open set $U$ such that 
$D \in \langle L \rangle$.
 
\end{itemize}
 
 \vspace{0.5cm}

Let us mention the following   
basic property of the symmetry operators though we shall not need it in the proof of the main result (Theroem 3).

 \vspace{0.5cm}

In the next lemma, we generalize the invariance which we studied in the case of the Euclidean Laplace operator in Example 1.

\vspace{0.5cm}
 
{\bf Lemma 1:} If $L$ is a symmetry operator of a differential operator $D,$ 
then $L$ preserves the kernel of $D,$ i.e.,
$L(\mbox{Ker} \, D) \subseteq \mbox{Ker} \, D.$

 \vspace{0.5cm}
 
{\it Proof.} If $Df=0$ for a smooth function $f$, we have 
$D(Lf) = L(Df)+[D,L]f = GDf = 0$ for a convenient operator $G.$
\hfill\(\Box\)

 \vspace{0.5cm}

{\bf Lemma 2}: Let $M$ be a manifold, $(U,\phi=(x^1,\ldots, x^n))$ be a map on this manifold, 
$D=\sum_{I \in K} f^I \partial_I \in DO(U)$ be a differential operator   
and $g(x^1,\ldots, x^n) = e^{\lambda x^1}$ for all $(x^1,\ldots, x^n) \in \phi(U).$ 
If for each $I \in K,$ $f^I$ does not  $\phi$-depend on the first variable $x^1$, then for every $\lambda \in \mathbb{C}$
the operator $D$ is $g$-invariant with respect to $(U,\phi).$

 \vspace{0.5cm}

{\it Proof.}   
Let $h \in C^{\infty}(U)$ be a function 
that does not $\phi$-depend on the first coordinate.
We have
\begin{equation*}
\begin{split}
    \left(\sum_{I \in K} f^I\partial_I\right)(e^{\lambda x^1}h) 
		&= \sum_{I=(i_1,\ldots, i_n) \in K}f^I 
\underbrace{\partial_{x^1}\ldots\partial_{x^1}}_{i_1 
\text{-times}}\partial_{I-(i_1,\:0,\ldots,\:0)}(e^{\lambda x^1}h) =\\
    &=\sum_{I=(i_1,\ldots, i_n) \in K} f^I \underbrace{\partial_{x^1}\ldots\partial_{x^1}}_{i_1 
\text{-times}}(e^{\lambda x^1} \partial_{I-(i_1,\:0,\ldots,\:0)}h)=\\
    &=\sum_{I=(i_1,\ldots, i_n) \in K}f^I \lambda^{i_1}e^{\lambda x^1}\partial_{I-(i_1,\:0,\ldots,\:0)}h=\\
&=e^{\lambda x^1}\sum_{I=(i_1,\ldots, i_n) \in K}\lambda^{i_1}f^I\partial_{(0,\ i_2,\ldots,\ i_n)}h.
\end{split}
\end{equation*}
Setting $D'=\sum_{I=(i_1,\ldots, i_n) \in K}\lambda^{i_1} f^I\partial_{(0,i_2,\ldots,\:i_n)},$ we see that 
$D$ 
is $e^{\lambda x^1}$-invariant with respect to $(U, \phi).$
\hfill\(\Box\)

\vspace{0.5cm}
 
{\bf Definition 4}: Let $(U, \phi)$ be a map on a manifold  $M.$ We say that a solution $f \in C^{\infty}(U)$ of $D f =0$ is 
 {\it $\phi$-separated} (or a $\phi$-separated solution) if
there exists an integer $1 \leq k < n$ and functions $g, h \in C^{\infty}(U)$ such 
that 
$f=gh$ where $g$ does not $\phi$-depend on the first $k$ variables and 
$h$ 
does not $\phi$-depend on the last $n-k$ variables.

 \vspace{0.5cm}

{\bf Remark:} Let us suppose that the dimension of the manifold $n>1.$ Then $f=0$ is always a $\phi$-separated solution. If $D$ does 
not contain the zero order term (in some and consequently in any coordinates), 
any constant function is a separated solution of $D.$

 Notice also that allowing $k=n$ in the above definition would lead
to the conclusion that any solution  $f$ is a $\phi$-separated solution since we might take $g$ to be an arbitrary non-zero constant function and set
$h=f/g.$  

If $n=1,$ we agree that no separable solution exists. 

\vspace{0.5cm}

 In the next theorem, we prove that the existence of a specific first order homogeneous simple 
symmetry operator implies the existence of a separated solution, which is non-constant. (The operator is actually a vector field commuting with 
$D.$)

 \vspace{0.5cm}
 
{\bf Theorem 3}: Let $M$ be a two dimensional manifold, $D \in DO(M)$ be a differential operator  and $L$ be a 
first order homogeneous simple symmetry operator of $D$ which does not generate $D.$
Then there is a point $m \in M$ and a map $(U,\psi)$ around $m$ such 
that there exists a non-constant  $\psi$-separated solution of $Df=0.$

 \vspace{0.5cm}

 {\it Proof.} Since $L$ is a first order homogeneous operator, it is  a non-zero vector field in particular.
Let us pick a point $m_0 \in M$ such that  $(Lf)(m_0)\neq 0$ 
for a smooth function $f$ on $M.$ 
By the theory of ordinary differential equations,   there exists  a map $(U_0, \phi=(x^1,x^2))$ 
around $m_0$ such that $L=\partial_{x^1}$ on $U_0.$ (The map can be defined using the local flow of $L$ and the flow  of a vector field $L'$
which may be taken   perpendicular to $L$  with respect to a  Riemannian metric defined around $m_0.$)

Let us   choose the connected component of $U_0$ containing 
$m_0,$  keeping denoted it by $U_0.$ Write the restriction of $D$ 
to $C^{\infty}(U_0)$ as $D=\sum_{|I|\leq d} f^I \partial_I$ (with respect 
to the map $\phi),$ where $f^I$ are smooth functions on $U_0$ and $d$ is a non-negative integer.
Since $L$ is a simple symmetry, we have $0 = [L,D] = [\partial_{x^1}, D] = \sum_{|I|\leq d} (\partial_{x^1}f^I) \partial_I$ because the coordinate vector fields commute with each other. In particular, $\partial_{x^1}f^I=0$ and thus, there exists an open set $U_1$ such that $f^I_{|U_1}$ does not $\phi$-depend on $x^1.$ 
By Lemma 2, $D$ is $g$-invariant for  $g(x^1,x^2)=e^{\lambda x^1}$ with respect to $(U_1, \phi).$

Because $D$ is $g$-invariant, we have the operator
\begin{equation*}
    D' = \sum_{j = 0}^{d_2} \left(\sum_{i = 0}^{d_1} \lambda^i f^{(i,j)} \right) \partial_{x^2}^j.
\end{equation*}
at our disposal. We take it in this form, that is derived in the proof of Lemma 2.
By assumption, $L=\partial_{x^1}$ does not generate 
$D=\sum_{i=0}^{d_1}\sum_{j=0}^{d_2}f^{(i,j)}\partial_{x^1}^i\partial_{x^2}^j.$ Thus 
$D$ has to contain at least one non-zero function coefficient $f^{(i, j)}$  in front of 
$\partial_{x^1}^{i}\partial_{x^2}^{j}$ for $0\leq i \leq d_1$ and $0 < j \leq d_2.$ 
Let us denote its indices by $i_0$ and $j_0$.

Let us consider a point $m_1\in U_1$ such that $f^{(i_0, j_0)}(m_1) \neq 0,$ and 
the polynomial $P(\lambda)=\sum_{i=0}^{d_1} f^{(i, j_0)}(m_1) \lambda^i,$ which has to be non-zero.
As a consequence of the continuity of $P$, there is an infinite set of complex numbers $\lambda$
for which this polynomial attains a non-zero value. 
Especially, it is possible to choose  $\lambda_0 \neq 0$ from this set.
Since $j_0\geq 1,$    $D'$ is a differential operator of order at least one. 
From now, we consider $D'$ only for $\lambda = \lambda_0.$

Since the chosen operator $D'$ contains derivatives only in the second variable and its function coefficients do not $\phi$-depend on $x^1$, $D'$ is a differential
operator in $x^2.$
By the basic theory of ordinary differential equations, $D' h = 0$ has a non-zero solution in an open neighbourhood 
$U \subseteq U_1$ of  $m_1$ since the degree of $D'$ is at least $1.$
We set $\psi=\phi_{|U},$ choose such a non-zero solution that does not $\psi$-depend on $x^1,$ and denote it by $\widetilde{h}.$ 
The function $f=g\widetilde{h}=e^{\lambda_0x^1}\widetilde{h}$  is non-constant. It is $\psi$-separated and it is a solution of $Df=0$ since
\begin{equation*}
    D(e^{\lambda_0 x^1}\widetilde{h}) = e^{\lambda_0 x^1} D'(\widetilde{h}) = 0.
\end{equation*}
It follows, that $(U, \psi=\phi_{|U})$ is a map around $m=m_1$ and $f=e^{\lambda_0 x^1}\widetilde{h}$ is a non-constant $\psi$-separated solution of 
$Df=0$ defined on $U.$ 
\hfill\(\Box\)

\bigskip

{\bf Remark:} 
\begin{itemize}
\item[1.] Instead of   the continuity of $P$ in the above proof, we could use the fact that 
non-zero polynomials in one variable have only a finite number of roots.
	
\item[2.] From the proof of Theorem 3, we see that if we can take even $j_0  > 1,$ there would be a non-constant 
solution $\widetilde{h}$ of $D'h=0,$ which does not $\phi$-depend on $x^1$. Thus
$f(x^1,x^2)=e^{\lambda x^1} \widetilde{h}(x^2)$ would be a  separated solution, whose first factor depends on $x^1$ only and the second one only on $x^2.$ 
More precisely, it is not true that the first factor does not $\phi$-depend on $x^1$ and it is not true that the second one does not 
$\phi$-depend on $x^2$.
\end{itemize}

\section{The Helmholtz operator}

As an application of Theorem 3, we derive a family (depending on a continuous parameter $\lambda)$ of classical  
separated solutions for the Helmholtz operator in two variables, i.e., for
\begin{align*}
    D=\partial_{x}^2 + \partial_{y}^2 + \omega^2, \mbox{ where } \omega \in (0,\infty)
\end{align*}
acting on smooth functions on $M =\mathbb{R}^2.$ Notice that the result is well known (see, e.g., \cite{Miller}).  
The operator $L = x\partial_y - y\partial_x$ is a first order homogeneous operator, which is a simple symmetry as the following computation shows

 \begin{align*}
   [x\partial_y - y\partial_x, D] =&  [x\partial_y - y\partial_x, \partial_x^2+\partial_y^2]\\
																	=&  (x\partial_y - y\partial_x)(\partial_{x}^2 + \partial_{y}^2) - (\partial_{x}^2 + \partial_{y}^2)(x\partial_y - y\partial_x)\\
																	=&  x\partial_y\partial_x^2+ x\partial_y^3 - y\partial_x^3 -y\partial_x\partial_y^2 \\
																	 &  -\partial_x(\partial_y+ x \partial_x \partial_y) +y\partial_x^3  - x\partial_y^3 + \partial_y (\partial_x  + y \partial_y\partial_x) \\
																	=&  x\partial_y \partial_x^2 + x\partial_y^3-y \partial_x^3  -y\partial_x\partial_y^2 - \partial_x\partial_y  - \partial_x\partial_y - \\
																	 &  - x\partial_x^2\partial_y + y\partial_y^3 -x\partial_y^3 + \partial_y \partial_x + \partial_y \partial_x + y\partial_y^2\partial_x =0.
 \end{align*} 
  We shall find the map $\phi$ and also a $\phi$-separated solution. 
For simplicity, let us consider only such points which do not lie on the non-positive part of the $x$-axis. (This limitation can be removed by 
translating  the found map by a constant vector.)  By the proof of Theorem 3, the symmetry operator $L$ has to be locally equal to 
a derivative with respect to a coordinate function of the map $\phi.$ It is easy to realize  that  $L = 
\partial_{\theta},$
where we use the polar patch $(r, \theta)$ in the form introduced in Example 1. (Obviously, $L$ does not determine the map $\phi$ uniquely.) In the considered polar patch, the Helmholtz operator  has the form $$D =\partial_{r}^2 + \frac{1}{r} \partial_r +\frac{1}{r^2} \partial_{\theta}^2  +  \omega^2.$$
Function coefficients of $D$ in the polar patch in front of the partial derivatives do not $\phi$-depend on the second variable, i.e., on $\theta.$
From this expression, we see that $L$ does not generate $D$ which contains a differentiation with respect to $r.$ According to Theorem 3, the Helmholtz operator has a non-constant separated solution  with respect to the map $\phi$ around any point in $M.$ 

  Using the construction from the proof of Theorem 3, we get
\begin{equation*}
(\partial_{r}^2 + \frac{1}{r} \partial_r +  \frac{1}{r^2} \partial_{\theta}^2 + \omega^2) (e^{\lambda \theta} h) =  e^{\lambda \theta} (\partial_{r}^2 +
\frac{1}{r} \partial_r +  \frac{\lambda^2}{r^2}+ \omega^2)h
\end{equation*}
where $h: r \mapsto h(r)$ does not $\phi$-depend on the last variable. 
We choose $D' =\partial_r^2 + \frac{1}{r} \partial_r 
+ \frac{\lambda^2}{r^2} + \omega^2.$ This operator is of degree $2$ independently on $\lambda$.  
By the theory of ordinary differential equations, we know that there exists not only a non-zero but also a non-constant solution of $D'h=0,$ i.e., of 
$(\partial_r^2 + \frac{1}{r}\partial_r +  \frac{\lambda^2}{r^2} + \omega^2) h = 0.$
Some of the solutions of this equation are  the well known {\it  Bessel functions of the first kind} $\widetilde{h}(r)=J_l(\omega r),$ $r>0,$ where $l= \pm \imath \lambda.$ (We omit a description of the whole fundamental system of solutions for $D'h=0,$ i.e., for the so called Bessel equation. See references below.)
All of these functions are non-constant, the function $J_0(\omega r)$ inclusively. (See, e.g., \cite{Lebedev},  \cite{MF} or \cite{Whitt}.)
Hence for every complex $\lambda \neq 0,$ the functions $(0,\infty) \times (-\pi,\pi) \ni (r,\theta) \mapsto f^{\pm}_{\lambda}(r, \theta) := e^{\lambda\theta}J_{\pm \imath \lambda}(\omega r)$ are $(r, \theta)$-separated solutions for the Helmholtz equation
$Df = 0$ on $M.$

\end{document}